\documentclass[10pt]{amsart}                          
\usepackage{epsfig}            
\usepackage{amssymb,amscd,bbm}     
\usepackage{amsthm,amsgen,amsmath,epic,psfrag}
\usepackage[dvips]{xypic}       

\xyoption{curve}                
\xyoption{rotate}               

\def\cup{\smallsmile}

\newcommand{\R}{{\mathbb R}}  
\newcommand{\Z}{{\mathbb Z}}  
\newcommand{\Q}{{\mathbb Q}}
\newcommand{\ka}{{\bf k}}
\newcommand{\I}{{\mathbb I}}       
\newcommand{\Lie}{\mathcal{L}\!{\it ie\/}}    
\newcommand{\Com}{\mathcal{C}\!{\it omm\/}}    

\newcommand{\C}{{\mathcal C}}

\newcommand{\dinv}{d^{-1}}

\newcommand{\dgaa}{{\text{\scshape dgaa}}}
\newcommand{\dgac}{{\text{\scshape dgac}}}
\newcommand{\dgc}{{\text{\scshape dgcc}}}
\newcommand{\dgl}{{\text{\scshape dgla}}}

\newcommand{\s}{ }

\newcommand{\lie}{\mathcal{L}\!{\it ie}}    

\theoremstyle{plain}                          
\newtheorem{theorem}{Theorem}[section]                          
\newtheorem{proposition}[theorem]{Proposition}                          
\newtheorem{lemma}[theorem]{Lemma}   
\newtheorem{corollary}[theorem]{Corollary}                          
    
\newtheorem{exercise}{Exercise}    
\newtheorem{problem}[exercise]{Problem}

\theoremstyle{definition}                          
\newtheorem{definition}[theorem]{Definition}  

\theoremstyle{remark}  
\newtheorem{example}[theorem]{Example}

\newcommand{\refT}[1]{Theorem~\ref{T:#1}}
\newcommand{\refC}[1]{Corollary~\ref{C:#1}}
\newcommand{\refP}[1]{Proposition~\ref{P:#1}}

\begin{document}

\title{Koszul duality in algebraic topology}
\author[D. Sinha]{Dev P. Sinha}
\address{Department of Mathematics\\
University of Oregon\\
Eugene, OR
97403}
\email{dps@math.uoregon.edu}

\maketitle

The most prevalent examples of Koszul duality of operads are the self-duality of the associative
operad and the duality between the Lie and commutative
operads.  At the level of algebras and coalgebras, the former duality was first noticed as such by
Moore,  as announced in his ICM talk at Nice \cite{Moor71}.
Thus this particular 
duality has typically been called Moore duality, and some prefer to call the
general phenomenon Koszul-Moore duality.  
The second duality at the level of algebras was realized in the
seminal work of Quillen on rational homotopy theory \cite{Quil69}.  Our aim in these notes based on 
our talk at the Luminy workshop on Operads in 2009 is to try
to provide some historical, topological context for these two classical algebraic dualities.  

We first review the original 
cobar and bar constructions used to study loop spaces and classifying spaces, emphasizing the
less-familiar geometry of the cobar construction.   Then, after some elementary topology,
we state duality between bar and cobar complexes in that setting.
Before explaining Quillen's work, 
we also share some other ideas - calculations of Cartan-Serre and Milnor-Moore and 
philosphy of Eckmann-Hilton - which may have influenced him.
After stating Quillen's duality, we share some recent work which relates these constructions
to geometry through Hopf invariants and in particular linking phenomena.

I would like to thank my collaborator Ben Walter.  This material is a union of standard material
which either I taught him or he taught me
along with new theorems which we have figured out together.

\section{Bar and cobar constructions}

\subsection{$\Omega X$ and the cobar construction}

Studying mapping spaces is one of the central tasks of topology, and
loop spaces are the simplest and most fundamental examples (unless one counts maps from finite 
sets, which  yield products).  We require a model for  loops where the loop sum
is associative exactly, not up to homotopy.  Thus, for us $\Omega X$ denotes the Moore loop
space which consists of pairs $f : \R \to X$ and a ``curfew'' $a > 0$ such that $f(x)$ is the basepoint
if $x \leq 0$ or if $x \geq a$.  Loop sum adds these curfews, which makes multiplication associative.

The cobar construction of Adams and Hilton \cite{AdHi56} was informed by the almost concurrent
work of  James \cite{Jame55} who studied $\Omega \Sigma X$, the loop space on the reduced suspension of $X$, namely $\Sigma X = X \times \I / (X \times 0 \cup * \times \I \cup 1 \times X)$.
There is a canonical inclusion of $J : X \hookrightarrow \Omega \Sigma X$ sending $x$ to 
$J(x)(t)$, the path which sends $t$ to 
the image in $\Sigma X$ of $(x, t)$.  Because $\Omega \Sigma X$ is a topological 
monoid, this map extends to a map from the free moinoid (with unit) on $X$  
to $\Omega \Sigma X$ which we call the James map $\hat{J}$.  
For example, the formal product $y * x * z$
goes to a loop with coordinates $(x,t)$ for $t \in [0,1]$ then $(y, t-1)$ for $t \in [1,2]$, then
$(z, t-2)$ for $t \in [2,3]$  -- see the figure below.

\begin{figure}\label{F1}
{
$$\includegraphics[width=3.5cm]{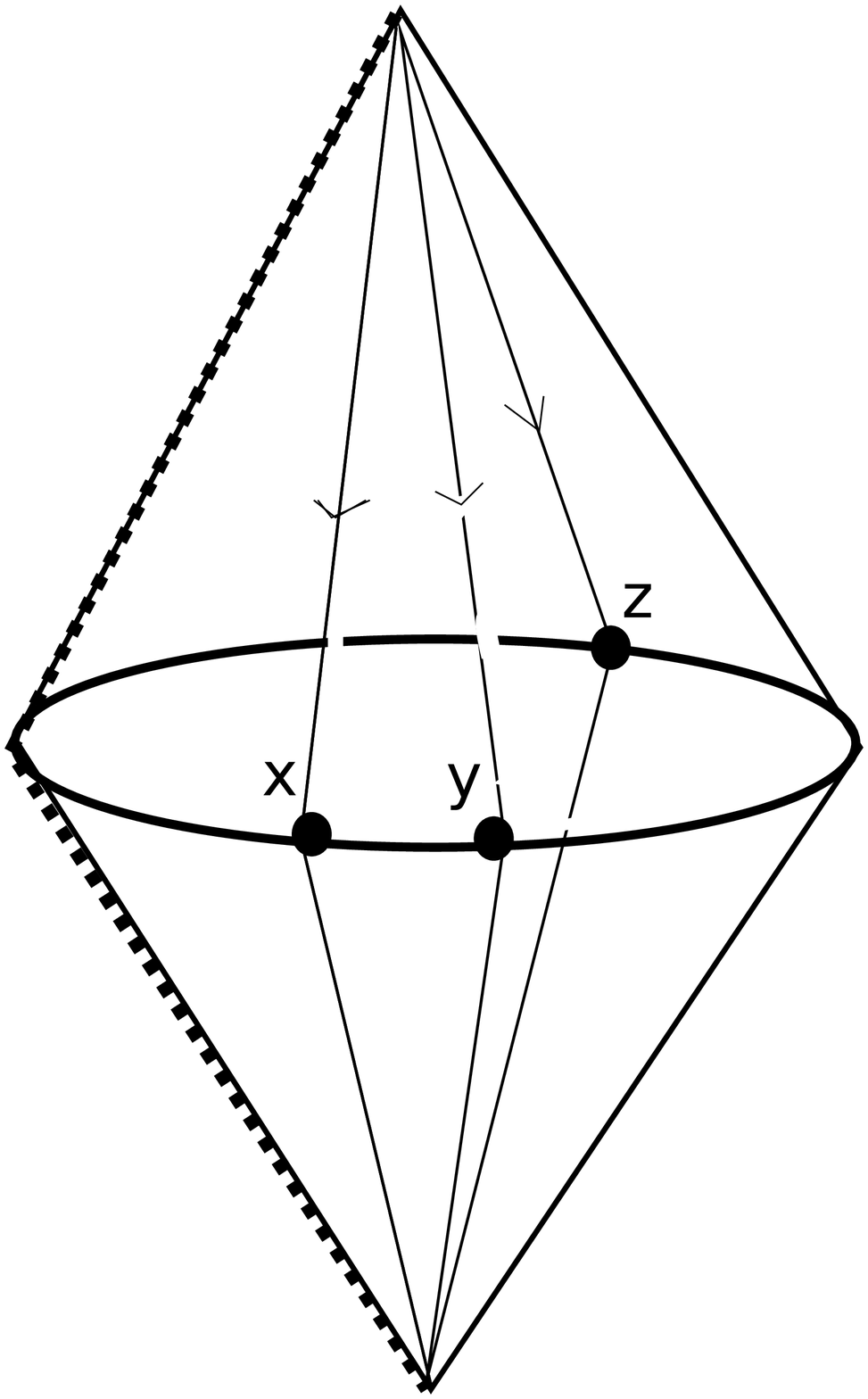}$$
}
An illustration of $\hat{J}$ of $y * x * z$ (traversing the path through $y$ first, etc).
\end{figure}

\begin{theorem} [James \cite{Jame55}]
The James map $\hat{J}$ from the free monoid on $X$ to $\Omega \Sigma X$ is a homotopy equivalence.
\end{theorem}

Recall that the homology of any space with an associative multiplication, or even a homotopy associative
multiplication, is an associative algebra.

\begin{corollary}
The homology of $\Omega \Sigma X$ with field coefficient is isomorphic as an algebra to the tensor (that is, free associative) algebra on the homology of $X$.
\end{corollary}

\begin{exercise}
Explicitly define the free topological monoid on a (well-based) topological space $X$.  Show
that its homology with field coefficients is isomorphic to the tensor algebra on the homology of $X$.
[Hint: make heavy use of the K\"unneth theorem.]
\end{exercise}

Interestingly, the corollary is typically proven in the course of proof of the theorem.  Details in a well-digested form are in Section~4.J in Hatcher's textbook \cite{Hatc02} 
or the survey paper of Carlsson and Milgram \cite{CaMi95}, whose treatment of the Adams-Hilton
construction heavily influences our  treatment below.  There are however more geometric proofs
which build on the fact that the space of paths in the cone on $X$ with endpoints in the image of
$X$ is homotopy equivalent to $X \times X$ through the projection  onto the endpoints.

\medskip

For the Adams-Hilton construction, we start with a simply-connected simplicial complex $\hat{X}$
and then contract the union of the $1$-skeleton along with enough of the two-skeleton
so that the quotient map $\hat{X} \to X$ is a homotopy equivalence.  Then $X$ is a CW-complex,
and its cellular chains are a quotient of the simplicial chains on $\hat{X}$.  By abuse of notation,
we denote these cellular chains by $C^{\Delta}_{*}(X)$.

Next, consider the cubical singular chain complex of the loop space $C^{\square}_{*}(\Omega X)$,
which is an associative differential graded algebra.  On generators, the product of 
$\sigma_{1} : \I^{n} \to \Omega X$ and $\sigma_{2} : \I^{m} \to \Omega X$ is the composite
$\I^{n + m} \cong \I^{n} \times \I^{m}  \overset{\sigma_{1} \times \sigma_{2}}{\to} 
\Omega X \times \Omega X \to \Omega X$.

The Adams-Hilton construction defines a map of associative algebras
from the free associative algebra on $C^{\Delta}_{*}(X)$
to $C^{\square}_{*}(\Omega X)$.  The first key observation is that any
choice of map $\gamma_{n} : \I^{n} \to \Delta^{n}$ defines a map  $AH_{\gamma_{n}} :
C^{\Delta}_{n}(X) \to C^{\square}_{n-1}(\Omega X)$.  Let $\chi_{\sigma}$ denote the characteristic
map $\Delta^{n} \to X$ of a simplex $\sigma$ of $X$.  Then $AH_{\gamma_{n}}(\sigma)$ is
basically given by the composite $\gamma_{n} \circ \chi_\sigma : \I^{n} \to X$.  From this composite 
we by adjointness (choosing say the last coordinate as the loop coordinate) 
produce a map $\I^{n-1} \to {\rm Map}(\I, X)$,
which then is identified with a generator of $C^{\square}_{n-1}(\Omega X)$ through viewing
${\rm Map}(\I, X)$ as Moore loops with curfew one.

The game is to define $\gamma_{n}$ appropriately so that  we can calculate boundaries, 
and more importantly so that the Adams-Hilton map
yields a quasi-isomorphism.  By abuse, we suppress $\gamma_{n}$ from notation
and write $AH_{\gamma_{n}}(\sigma)$ as $|\sigma|$. 
For the first case when $n=2$, a good way to choose $\gamma_{2}$ is to
to have $\gamma_{2} : \I^{2} \to \Delta^{2}$ send the boundary of $\I^{2}$ to that of $\Delta^{2}$.
In any way this is done, we would have that $d|\sigma_{2}| = |d \sigma_{2}| = 0$, since the one-skeleton
of $X$ is has been identified to a point.  Looking forward, it is much better to choose $\gamma_{2}$
to be a ``degree one'' map $\I^{2} \to \Delta^{2}$ which when we consider the adjoint 
$\hat{\gamma}_{2} : \I \to {\rm Map}(\I, \Delta^{2})$ interpolates 
between the direct path from vertex $0$ to vertex
$2$ of $\Delta^{2}$ along the edge between them and the ``long'' path from $0$ to $2$ which first
traverses the $0$-$1$-edge and then the $1$-$2$ edge.  When composed with the 
characteristic map into $X$ these edge paths will yield constant loops, but the choice of the paths
in between is important.

\begin{figure}\label{F2}
$$\includegraphics[width=6cm]{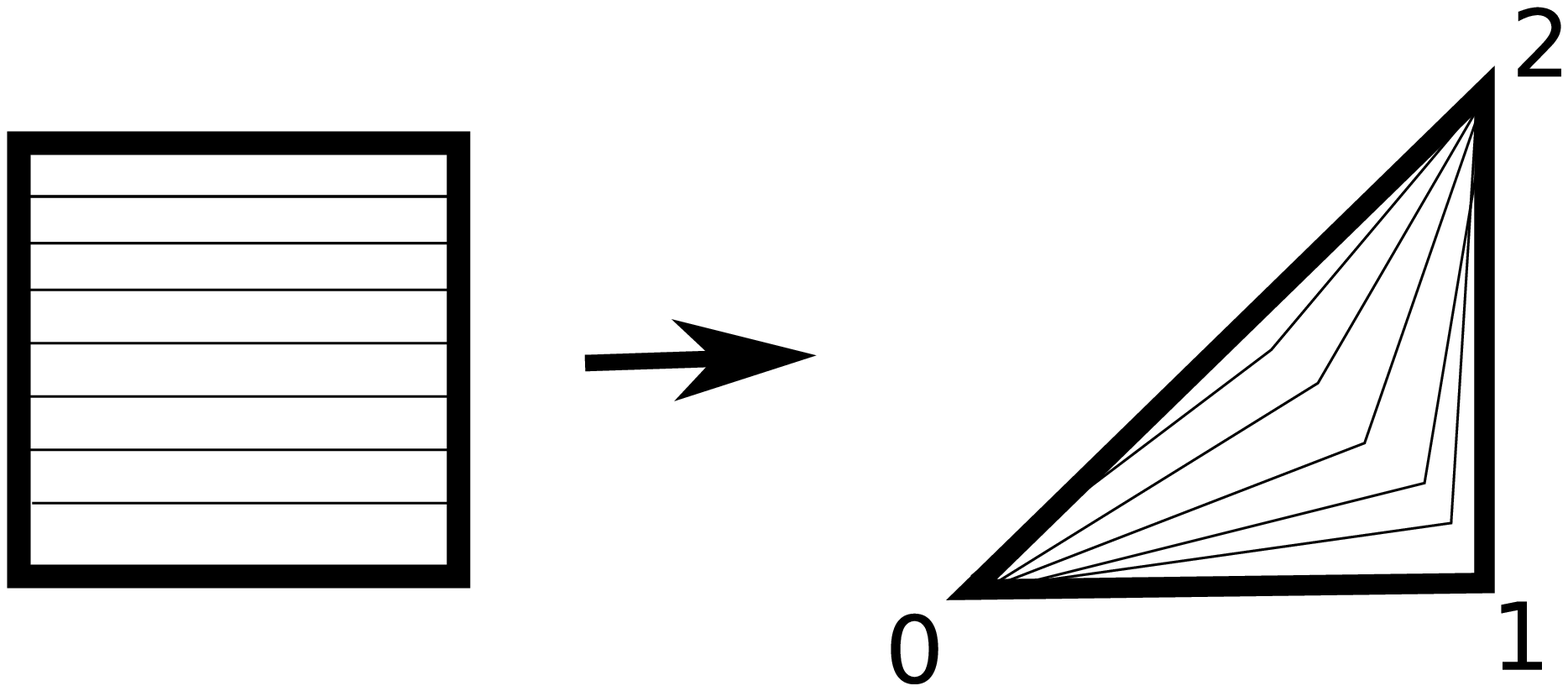}$$
One possible choice for $\gamma_{2}$.
\end{figure}

At the next stage, building on some such choice of degree one $\gamma_{2}$, we can define
a $\gamma_{3}$ such that $d|\sigma_{3}| = |d \sigma_{3}|$.  There are four faces
of $\I^{2}$ and four faces of $\Delta^{3}$, and this equality identifies those faces.  For example,
one face of $\I^{2}$  will be mapped by $\hat{\gamma}_{3} : \I^{2} \to {\rm Map}(\I, \Delta^{3})$ to 
paths from the $0$ vertex to the vertex $3$ (reminder: such a path in $\hat{X}$ will project to a loop 
in $X$) which first go to the vertex $1$ directly along the $0$-$1$ edge
and then go to $3$ along paths
compatible with the choice made of $\gamma_{2}$.   On another face of $\I^{2}$,
paths go only along the $0$-$2$-$3$ face of $\Delta^{3}$, again compatibly with $\gamma_{2}$,
and so on.

At $n= 4$ the construction there is a surprise.  Assume $\gamma_{3}$ has been defined,
and start defining $\hat{\gamma_{4}} : \I^{3} \to \Delta^{4}$ 
by setting its restriction to various faces as before,
for example sending one face of $\I^{3}$ to paths  on the $0$-$2$-$3$-$4$ face of 
$\Delta^{4}$.  But, there are six faces of $\I^{3}$ and only five faces of $\Delta^{4}$!  What is the natural
last term?  When one does the geometry carefully, one see that on the last face of $\I^{3}$ should
map to paths from $0$ to $4$ which first go along the $0$-$1$-$2$ face and then along the
$2$-$3$-$4$ face.  These two faces appear in the standard definition of the coproduct on simplicial
chains dual to cup product.  Moreover, such composites are given by the product in $C^{\square}_{*}(\Omega X)$.  That is, we can construct $\gamma_{4}$ such that
$$ d|\sigma_{4}| = |d \sigma_{4}| + |\alpha_{2}| * |\beta_{2}|, $$
where the coproduct of $\sigma_{4}$ is $\alpha_{2} \otimes \beta_{2}$ plus terms
in bidegrees $(1,3)$, $(0,4)$, etc.  (These other terms in the coproduct  yield trivial chains with the $1$-skeleton of $X$ collapsed.)  

In general, let us denote products of Adams-Hilton chains by $|\sigma| * |\tau| = |\sigma|\tau|$.
Let ${\rm Cobar}(C^{\Delta}_{*}(X))$ denote the sub-algebra of $C^{\square}_{*}(\Omega X)$ generated
by the Adams-Hilton chains (in positive degrees).

\begin{theorem}[Adams-Hilton]\label{T:AH}
There are degree-one choices for the maps $\gamma_{n}$ such that the boundary on 
${\rm Cobar}(C^{\Delta}_{*}(X))$  is the cofree extension of the map with 
$$ d |\sigma| = \pm |d \sigma| + \sum_{\bar{\Delta} \sigma = \sum \alpha_{i} \otimes \beta_{i}}
\pm  |\alpha_{i}|\beta_{i}|.$$
Here $\bar{\Delta}$ denotes the reduced cup coproduct including terms of only positive bidgrees.

The inclusion of any such ${\rm Cobar}\left(C^{\Delta}_{*}(X)\right)$  in $C^{\square}_{*}(\Omega X)$ is a 
quasi-ismorphism of differential graded associative algebras.  
\end{theorem}

\begin{exercise}
Try to write down $\gamma_{n}$ for $n \leq 4$ as an explicit piecewise-linear map.
\end{exercise}

\begin{exercise}
The cobar construction is defined for any differential graded coalgebra.  Compute it for the coalgebra
given by the homology of $\C P^{\infty}$.
\end{exercise}

\begin{exercise}
Deduce the James theorem from the Adams-Hilton theorem.
\end{exercise}

The algebraic 
cobar construction (often denoted $\Omega$ but not at the moment because of
potential confusion) has become part of the standard toolkit for algebraic topologists, and there
are more algebraic approaches which can yield similar theorems.  A more geometric approach to
the topology of iterated loop spaces was extended by Milgram who studied $\Omega^{n} \Sigma^{n} X$
in \cite{Milg66} (see also \cite{CaMi95}).  
But the geometry and formalism of PROPs and operads, in particular the elegance of the little
disks construction of Boardman and Vogt \cite{BoVo73}, became more popular than this intricate
geometry.  Perhaps there could be something gained by revisiting these ideas.

\subsection{Classifying spaces and the bar construction}

We will be more brief about the bar construction, whose topology is better known.
The topological bar construction provides a model for the classifying space $BG$,
which when $G$ is discrete is just the Eilenberg-MacClane space $K(G, 1)$.

Topologists are often ambiguous and refer to any quotient of
a contractible space with free $G$-action as its classifying space
$BG$.  We resolve this issue by only saying that 
a space is homotopy equivalent (rather than equal to) $BG$. As we remark below, there is
a choice which is useful to call $BG$ unambiguously.

\begin{example}
\begin{itemize}
\item $B\Z \simeq S^{1}$.
\item $B\Z/2 \simeq \R P^{\infty}$.
\item $B\Z/n \simeq S^{\infty}/ (\Z/n)$, called an infinite Lens space.
\item If $G = \pi_{1}(S)$ where $S$ is a surface of positive genus, then $BG \simeq S$.
\end{itemize}
\end{example}

\begin{theorem}\label{T:BG}
If $G$ is discrete, then $BG$ is homotopy equivalent to a simplicial complex whose $n$-simplices
are in one-to-one correspondence with $n$-tuples of elements of $G$, which we denote
$|g_{1}|g_{2}| \cdots |g_{n}|$.  The $(n+1)$ faces of an $n$-simplex are given by 
$$ d_{i}(|g_{1}| \cdots |g_{n}|) = 
\begin{cases}
|g_{2}| \cdots | g_{n}|    & i = 0\\
|g_{1}| \cdots |g_{i} g_{i+1} |  \cdots |g_{n}| & 0 < i < n \\
|g_{1}| \cdots | g_{n-1} |  & i = n.
\end{cases}
$$
\end{theorem}

To prove this, one constructs $EG$ in a similar fashion.

\begin{corollary}\label{C:BG}
The homology of $BG$ is given by the homology of the algebraic bar construction applied
to the group ring $\ka[G]$, an associative algebra.
\end{corollary}

\begin{exercise}
Do the simple unraveling of definitions to check that this corollary follows.
\end{exercise}

We obtain a better model if we quotient by identifying each $n$-simplex  of the
form $|g_{1}| \cdots |e | g_{i+1} |
\cdots | g_{n} |$ with the $(n-1)$-simplex $|g_{1}| \cdots |g_{i-1} | g_{i+1} | \cdots | g_{n} |$ through the
appropriate standard projection of $\Delta^{n} \to \Delta^{n-1}$.  
The following exercise is a must for any topology student.

\begin{exercise}
Show that this reduced construction for $\Z/2$ is homeomorphic to  $\R P^{\infty}$.  \end{exercise}

Thus  $\R P^{\infty}$ has $\Z/2$ as its DNA, so to speak.
\refT{BG} is true in greater generality in particular when $G$ has a topology (with some
mild assumptions) which gets incorporated in the topology on $BG$,  or when
$G$ is just a monoid. Indeed, this construction is a special case of the
nerve of a category.

\subsection{Relating the bar and cobar constructions}

We defined the homotopy type of $BG$ through the fiber sequence
$$ G \subset EG \to BG.$$
Let $PX$ denote the path space on $X$, which is contractible, and let $ev$ denote the map which 
sends a path $\gamma$ to $\gamma(1) \in X$.  Then the sequence
$$ \Omega X \to PX \overset{ev}{\to} X$$
is a fibration.
Consider as well the map $P EG \to BG$ defined by evaluation composed by the quotient.  This map is equivalent
to both the projection $EG \to BG$ and the evaluation $PBG \to BG$, which are thus equivalent to each other.  We deduce
that their fibers are equivalent, so that $\Omega BG \simeq G$.
Similarly, if $X$ is connected then $B \Omega X \simeq X$ (the content of this statement
depends on the definition of classifying space for $\Omega X$; 
some say its classifying space is $X$ by definition).

These homotopy equivalences are reflected in the following algebra, which is now viewed as a consequence of Koszul
duality of the associative operad.  Recall that the cobar construction was defined in 
terms of a free associative algebra (and indeed computed the homology of $\Omega X$ as an
algebra).  We can view the bar construction as based on the free coassociative coalgebra generated by $\ka[G]$, with the coproduct  defined by breaking bar expressions in two and differential defined using the product of $G$.  

\begin{theorem}
The bar construction $B$ and the cobar construction $\Omega$ define
an adjoint pair of functors between differential graded associative algebras $\dgaa$ and
differential graded associative coalgebras $\dgac$. 
\begin{equation*}\label{bigdi}
\begin{aligned}
\xymatrix@=20pt{
\dgac \ar@<.5ex>[r]^{\Omega} & 
  \dgaa \ar@<.5ex>[l]^{B}  
}
\end{aligned}
\end{equation*}
 Moreover, there are natural
transformations $\Omega B A \to A$ and $B \Omega C \to C$ which if are quasi-isomorphisms
if $A$ is positively graded and if $C$ is $1$-connected respectively.
\end{theorem}

This theorem was announced by Moore \cite{Moor71}, so it has historically been referred to as Moore duality.  
In topology, this equivalence reflects the bijection between homotopy classes of 
monoid maps from some $M$ to $\Omega X$ and homotopy classes of maps
from $BM$ to $X$.

Not only is it the first example of adjoint functors giving equivalences between categories of
algebras and coalgebras over an operad and its Koszul dual, but it played a central role in Priddy's
definition of Koszul quadratic algebras \cite{Prid70}.  A graded augmented algebra $A$ can be given a zero differential.  Over a field $\ka$ and with finiteness degree-wise, 
the homology of the bar complex of $A$ is the linear dual of ${\rm Ext}_{A}(\ka, \ka)$, compatible with
their coalgebra and algebra structures.  (In the case
of $A = \ka[G]$, this is reflected by \refC{BG} and the fact that the cohomology of $BG$ is coincides
with ${\rm Ext}_{\ka[G]}(\ka, \ka)$.)   If $A$ is a Koszul algebra, then we can replace the bar complex 
with a much smaller resolution, which leads to an explicit presentation of this ${\rm Ext}$-algebra.
Moreover, the theory applies to this ${\rm Ext}$-algebra as well and replaces the cumbersome
quasi-isomophism of $A \simeq \Omega B A $ with an isomorphism 
$A \cong {\rm Ext}_{{\rm Ext}_{A}(\ka, \ka)}(\ka, \ka)$.

\section{Other ideas in the air}

Following up on his thesis, Serre along with Cartan considered the rational homotopy groups
of a simply connected space.  When shifted down, as best done by considering the
homotopy groups of $\Omega X$, those groups form a graded Lie algebra.
Typically the Hurewicz homomorphism from homotopy to homology captures little information.
But rationally for loop spaces, this map gives a clear picture.  Building on calculations of Cartan
and Serre \cite{CaSe52},  Milnor and Moore in \cite{MiMo65} prove the following.

\begin{theorem}\label{T:mimo}
If $X$ is simply connected,
the Hurewicz map $\pi_{*}(\Omega X)\otimes \Q  \to H_{*}(\Omega X; \Q)$ is an injection,
mapping the rational
homotopy Lie algebra of $X$  to the primitives in the Hopf algebra $H_{*}(\Omega X; \Q)$.
\end{theorem}

\medskip

Another influential idea at that time was Eckmann-Hilton ``Duality,'' which draws attention
to parallel structures in cohomology and homotopy.    See the table below.

\begin{table}
\begin{center}
\begin{tabular}{c|c}
Cohomology & Homotopy \\
\hline \\
L.E.S of a cofibration $A \hookrightarrow X \to X/A$ &
L.E.S of a fibration $F \to E \to B$ \\  \\
Spheres and Moore spaces & Eilenberg-MacClane spaces \\ \\
Suspension / desuspension  &  Loop space  / classifying space \\ \\
CW structures  & Postnikov tower \\ \\
Graded commutative ring structure & Graded Lie algebra structure \\ \\
co-$H$-space (comonoid) &  $H$-space (monoid) \\  \\
pushout square / homotopy colimit & pull-back square  / homotopy limit \\  \\
Steenrod algebra & Stable homotopy groups of spheres \\  \\
Leray-Serre spectral sequence  & Blakers-Massey theorems  
\end{tabular}
\end{center}
\end{table}

This duality is more of a philosophy than a theory.  There are no theorems of the form
``Given a true statement about homotopy groups, there is a true statement about cohomology
groups obtained by...''  or ``Given a space $X$ there is a dual space $\hat{X}$ whose cohomology
groups are the homotopy groups of $X$ and...''  Nonetheless, the duality can point to interesting
directions of study.  For example, looking at our table one notices a significant difference between
CW structures, which are not canonical in any sense, and the Postnikov tower, which is.  This
leads to finding the homology decomposition of a space (see Chapter~4.H of \cite{Hatc02} for
a basic treatment).

\section{Quillen functors and rational homotopy theory}

Quillen, influenced by Kan, took the step in \cite{Quil69} of proving theorems not about homotopy groups
but about all of homotopy theory.  He must have taken Theorems \ref{T:mimo} and \ref{T:AH}
as an important starting point.  Indeed, if the rational homology of the cobar construction computes 
the homology of the loop space, and one is to then take primitives to get rational homotopy groups, 
why not take primitives first at the level of the cobar complex itself (see exercise below)?   The great
advantage is that in the cobar construction one is considering the {\em free} associative algebra, whose
primitives are known to be the free Lie algebra, so one can just use the free Lie algebra 
functor  as a starting point.  
Quillen was also aware of Chevalley-Eilenberg cohomology of Lie algebras 
\cite{ChEi48}, and probably knew of some cases in which applying this functor to the rational homotopy
Lie algebra of a space recovered its cohomology (an easy case being wedges of spheres,
whose rational homotopy Lie algebra is free).    Once again, a refinement is needed, going from
applying a functor at the level of algebras (in the previous case primitives, in the current
case Lie algebra cohomology) to 
applying it at the level of chain complexes.  Quillen's adaption of the Chevalley-Eilenberg
construction now bears his name as well.  

Quillen put these two constructions together in the following theorem.

\begin{theorem}\label{T:main}
The  Lie algebraic cobar construction $\Omega_{\Lie}$ and a commutative coalgebraic
bar construction $B_{\Com}$, which generalizes the Chevalley-Eilenberg construction,
form an adjoint pair of functors 
\begin{equation*}
\begin{aligned}
\xymatrix@=20pt{
\dgc \ar@<.5ex>[r]^{\Omega_{\Lie}} & 
  \dgl \ar@<.5ex>[l]^{B_{\Com}}  
}
\end{aligned}
\end{equation*}
Here $\dgc$ are $1$-connected differential graded cocommutative coalgebras and 
$\dgl$ are connected differential graded Lie algebras.
These functors preserve all notions relevant to homotopy theory (fibrations, cofibrations, weak equivalences).

Any simply-connected 
space $X$ has functorial models $C_{X}$ and $L_{X}$ in $\dgc$ and $\dgl$ respectively
such that the homology of $C_{X}$ is the rational homology coalgebra of $X$ and the homology
of $L_{X}$ is the rational homotopy Lie algebra of $X$.
\end{theorem}

In current language, we would say that $\Omega_{\Lie}$ and $B_{\Com}$ form a 
Quillen adjoint pair of functors on the model categories $\dgc$ and $\dgl$, reflecting
the Koszul duality of the operads $\Lie$ and $\Com$.
This theorem gives a precise manifestation of Eckmann-Hilton duality, through the fact that
these functors preserve model structures along with the symmetries of the model  structure
axioms.
What complicates \cite{Quil69} significantly is that there is, to this day, no simple way
to construct a commutative cochain algebra of a space and thus easily  land in this picture.
Quillen has to walk for forty days through the desert, producing a long chain of functors in order
to produce $L_{X}$ and $C_{X}$.  That difficulty led Sullivan to find a simple way to produce
commutativity, not on chains but on cochains.  Additionally, 
instead of using bar or cobar constructions  Sullivan  studied 
cofibrant replacements with some additional smallness property, the famous minimal models
of  \cite{Sull77}.

\begin{exercise}
Check directly in some cases
that the primitives of differential graded Hopf algebra form a split sub-complex,
so that the primitives of the homology of $C_{\bullet}$ is isomorphic to the homology of the complex
obtained by taking the primitives of $C_{\bullet}$.
\end{exercise}

\begin{exercise}
Compute the Chevalley-Eilenberg cohomology of the graded Lie algebra with three
generators $x,y,z$ in degree three with the only relation being $[x,y] = [y,z]$.
\end{exercise}

\section{Koszul duality and Hopf invariants}

We have recently found \cite{SiWa08}
that Koszul duality and Quillen functors allow one to give a
definitive treatment of rational homotopy functionals through Hopf invariants. 
The basic idea can be seen as using the bar complex to understand a map
$f : S^{n} \to X$ by first passing to $\Omega f : \Omega S^{n}  \to \Omega X$ and
then evaluating cohomology classes of $\Omega X$ on the image of the fundamental
class of $\Omega S^{n}$.  By \refT{mimo}, such invariants are complete.

We must pause to make a choice in notation.  If one is studying the cohomology of $\Omega X$
using the cochains on $X$, one could either denote the construction you use by $\Omega$ to 
reflect topology or $B$ to denote a bar construction which is applied to an algebra (rather
than a cobar construction applied to a coalgebra).  The algebraists seem to have won this notational
conflict, so we consider $B C^{*}(X)$, the bar construction on the cochains of $X$ with their
associative cup product.  We let $H_{B}^{*}(X)$ denote the homology of $B C^{*}(X)$.
Define the weight of a generator of a bar complex to be the number of elements appearing.

\begin{lemma}\label{L:key}
$H^{n-1}_{B}(S^n)$ is rank one, generated by an element of weight one corresponding
to a  generator of $H^n(S^n)$. 
\end{lemma}

\begin{exercise}
Prove this.  Hint: you'll need the K\"unneth theorem to put yourself in a position to do
some ``weight reduction,'' as we use below.
\end{exercise}

\begin{definition}\label{D:int_B}
Let $\gamma \in B^{n-1}(C^{*} S^n)$ be a cocycle.  Define $\tau(\gamma) \simeq \gamma$ 
to be a choice of weight one cocycle to which $\gamma$ is cohomologous.

Define $\int_{B(S^n)}$ to be the map from cocyles in $B^{n-1}(C^{*}S^n)$ to $\mathbb{Z}$ given by
$\int_{B(S^n)} \gamma = \int_{S^n} \tau(\gamma)$, where $\int_{S^n}$ denotes evaluation
on the fundamental class of $S^n$.

Define $\eta(\gamma)$, the Hopf invariant associated to $\gamma$ by $\eta(\gamma) (f)
= \int_{B(S^{n})} f^{*}\gamma$.
\end{definition}

The choice of Hopf cochain is not unique, but the corresponding Hopf invariants are.  It is
immediate that the Hopf invariants are functorial.   Moreover, note that the definitions
hold with any ring cofficients.  Topologically we have the following
interpretation.

\begin{proposition}\label{P:spacelevel}
$\eta(\gamma)(f)$ coincides with the evaluation of the cohomology class given
by $\gamma$ in $H^{n-1}(\Omega X)$ on the image under $\Omega f$ of the
fundamental class in $H_{n-1}(\Omega S^{n})$.
\end{proposition}

\subsection{Examples}

\begin{example}
A cocycle of weight one in $B(C^{*}X)$ is just a closed cochain on $X$, which may be pulled back and 
immediately evaluated.  
Decomposable elements of weight one in $B(X)$ are null-homologous, 
consistent with the fact that products evaluate trivially on the Hurewicz homomorphism.
\end{example}

\begin{example}\label{e2}
Let $\omega_{1}$ and $\omega_{2}$ be generating $2$-cocycles on $S^{2}$ and $f : S^{3} \to S^{2}$.  
Then $\gamma =  -|\s\omega_{1} | \s\omega_{2}|$ is a cocycle in $B(C^{*}S^{2})$ which $f$ pulls back
to $- |\s f^{*} \omega_{1} | \s f^{*} \omega_{2} |$, 
a weight two cocycle of total degree two on $S^{3}$.
Because $f^{*} \omega_{1}$ is closed and of degree two on $S^{3}$, it is exact. 
Let $\dinv f^{*} \omega_{1}$ be a choice of a cobounding cochain.  Then 
$$d \left(|\s\dinv f^{*} \omega_{1} | \s f^{*} \omega_{2}| \right) = 
\s |f^{*} \omega_{1} |  \s f^{*} \omega_{2}| \ +\ 
\s |\dinv f^{*} \omega_{1} \cup f^{*} \omega_{2}|.$$
Thus $f^*\gamma$ is homologous to $\s|\dinv f^{*} \omega_{1} \cup f^{*} \omega_{2}|$,
and the corresponding Hopf invariant is $ \int_{S^{3}} {\dinv f^{*} \omega_{1} \cup f^{*} \omega_{2}}$,
which when choosing $\omega_{1} = \omega_{2}$ 
is  the classical formula for Hopf invariant given by Whitehead \cite{Whit47}.
It is a direct translation of the linking number definition of Hopf invariant into the
language of cochains.  Understanding this formula from the point of view of the bar
construction has, to our knowledge, only come over fifty years since all of these concepts were
introduced.
\end{example}

\begin{figure}\label{F3}
\psfrag{w}{$\omega_1$}
\psfrag{v}{$\omega_{2}$}
\psfrag{b}{$f^{*}\omega_{2}$}
\psfrag{a}{$f^{*}\omega_{1}$}
\psfrag{c}{$d^{-1}(f^{*}\omega_{1})$}
\psfrag{d}{$d^{-1}(f^{*}\omega_{1}) \cup f^{*}\omega_{2}$}
$$\includegraphics[width=12cm]{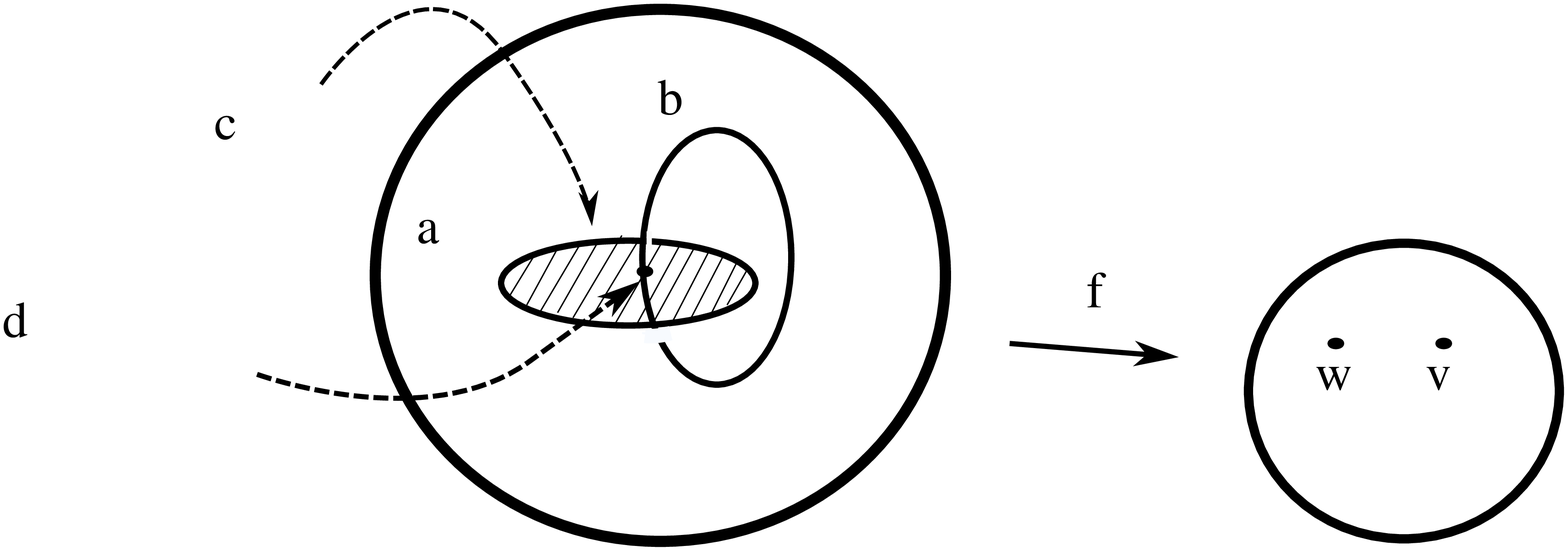}$$
Whitehead's integral, viewed through intersections of supports of cochains.
\end{figure}

\begin{example}\label{E:arbwt1}
For an arbitrary $X$ and cochains $x_{i}, y_{i}$ and $\theta$ on $X$ with 
$d x_{i}  = d y_{i} = 0$ and $d \theta = \sum (-1)^{|x_i|} x_{i} \cup y_{i}$, the cochain
$\gamma = \sum |\s x_{i} | \s y_{i} |+ |\s\theta| \in B(C^{*}X)$ is closed. 
The possible formulae for the Hopf invariant are all of the form  
$$ \int_{S^{n}} \left( f^{*} \theta - 
   \sum \left((-1)^{|x_i|} t_{i} \cdot \dinv f^{*} x_{i} \cup f^{*} y_{i} 
         + (1-t_{i}) \cdot f^{*} x_{i} \cup \dinv f^{*} y_{i}\right) \right),$$
for some real numbers $t_{i}$.
This generalizes  formulae given in \cite{Sull77}, \cite{GrMo81}
and \cite{Haef78}.

By choosing $t = \frac{1}{2}$ we see that reversing the order to consider $\sum |\s y_{i} | \s x_{i}|$
will yield the same Hopf invariant, up to sign.  Thus $\sum| x_{i} | y_{i} | \mp | y_{i} | x_{i}|$ yields
a zero Hopf invariant. There are many Hopf invariants of the classical bar construction
which are zero, a defect remedied by using the Lie coalgebra cobar construction.
\end{example}

\begin{exercise}
Suppose $x$ and $y$ are cochains supported on codimension two submanifolds $X$ and
$Y$ of $W$ and $\theta$ satisfies $d\theta = x \cup y$ and 
is supported on a codimension three submanifold which cobounds $X \cap Y$.
Draw pictures of how the Hopf invariant associated to $|x|y| \mp |\theta|$ 
evaluates some map $f : S^{3}  \to X$.  Moreover, draw pictures of what can happen
in $S^{3} \times \I$ if one has a homotopy between $f$ and $g$.  [Hint:  Start with the picture
in the figure, but then draw in the preimage of the support of $\theta$;  then, think about
what can happen with the preimage of $X \cap Y$ through a homotopy.]
\end{exercise}

One can do similar calculations in higher weight, and interpret them all when one chooses
cochains supported on submanifolds in terms of linking behavior of the preimages of those
submanifolds.  See \cite{SiWa08}.

\subsection{The cokernel and kernel of the Hopf invariant map}

Our Hopf invariant construction defines  a homomorphism
$\eta: H_{*} (B (C^{*}(X; \Z))) \to {\rm Hom}(\pi_{*} X, \Z)$.  
It follows from \refP{spacelevel} and \refT{mimo}
that this map is surjective when tensored with the rational
numbers, and thus is full rank.

\begin{problem}
Compute the cokernel of $\eta$.  By Adams' celebrated result \cite{Adam60}, this cokernel is trivial 
for $X$ an odd sphere and for $S^{2}$, $S^{4}$ and $S^{8}$, and it is $\Z/2$ for other even
spheres.
\end{problem}

The proofs in \cite{SiWa08} show that one might be able to 
directly understand the relation of this cokernel to lack of commutativity of cup product.
Though this cokernel is clearly a very subtle homotopy invariant,  we do not see any applications
of its calculation.  

Also, $\eta$ has a very large kernel, 
explained from the operadic viewpoint as the fact that we are taking the wrong 
bar construction.  The rational PL cochains on a simplicial set are commutative,
so we should be taking a bar construction over the Koszul dual cooperad, namely the 
Lie cooperad, rather than associative cooperad. The homology of 
such a bar construction $B_{\lie}$ is known as Harrison homology.  Using a graphical
model for the Lie cooperad developed in \cite{SiWa06} which makes calculations
possible, we prove the following.

\begin{theorem}\cite{SiWa08} \label{T:complete}
There is a Hopf invariant map $\eta^{\lie}$ which factors the map $\eta$ such that 
$\eta^{\lie} : H^{*}_{B_{\lie}}(X) \to {\rm Hom}( \pi_{*}(X),\; \Q)$ is an isomorphism
of Lie coalgebras.
\end{theorem}

It is almost immediate that similar Hopf invariants can be used to concretely 
realize similar isomorphisms arising for Koszul pairs in general.  


To summarize, in homology theory it has been helpful to have geometry attached
not only to homology but cohomology.  In particular, homology classes are often represented
by closed submanifolds and cohomology classes are represented by either forms or proper submanifolds.
The geometry of homotopy groups arising from their definition is almost too simple.
To reflect on the geometry of \refT{mimo}, we notice that the Lie algebra generators
of $\pi_{*}(X) \otimes \Q$ have non-trivial Hurewicz image, as noticed by Cartan-Serre.
That is, the rational homotopy groups of $X$ are spanned by Whitehead products of
spherical homology classes.  Our work on Hopf invariants 
shows that the geometry of homotopy functionals
is given by linking invariants, as perfectly governed by the Lie cooperad, completing the
geometric understanding of these basic functors in the rational setting.  We hope these
ideas can be extended to the non-simply connected setting, and perhaps - at least in part - 
in characteristic $p$.

\bibliographystyle{amsplain}
\bibliography{references}

\end{document}